\documentclass{article}
\usepackage{amsmath, amsthm, amsfonts}

\theoremstyle{definition}

\theoremstyle{remark}

\input{tcilatex}

\begin{document}

\title{An elliptic partial differential equation and its application}
\author{Dragos-Patru Covei \\
{\small The Bucharest University of Economic Studies, Bucharest, Romania. }\\
{\small E-mail address: dragos.covei@csie.ase.ro} \and Traian A. Pirvu \\
{\small McMaster University, Hamilton, Canada. }\\
{\small E-mail address: tpirvu@math.mcmaster.ca}}
\maketitle

\abstract{This paper deals with the following elliptic equation\begin{equation*}
-2\left\vert \sigma \right\vert ^{2}\Delta z+\left\vert \nabla z\right\vert
^{2}+4\alpha z=4\left\vert x\right\vert ^{2}\text{ for }x\in \mathbb{R}^{N}\text{, }(N\geq 1),
\end{equation*}where  $\alpha >0,$ and $\sigma>0$ are some
real constants. The solution method is based on the sub- and super-solution
method. The case $N>1$ seemed not considered before. This equation models a
stochastic production planning problem.}

\section{Introduction}

In this article, we look for positive solutions of the following partial
differential equation

\begin{equation}
-2\sigma^{2}\Delta z+\left\vert \nabla z\right\vert ^{2}+4\alpha
z=4\left\vert x\right\vert ^{2}\text{, }x\in \mathbb{R}^{N}\text{.}
\label{tpes}
\end{equation}%
Here $N\geq 1$ is the space dimension, $\left\vert \circ \right\vert $ is
the Euclidean norm of $\mathbb{R}^{N}$, $\alpha >0,$ and $\sigma>0$ are some
real constants.

This equation has received much attention in the last decades since it is
related with several models that arises in production planning problems; for
more on this see the papers of Akella and Kumar \cite{AK}, Alvarez \cite{A},
Bensoussan, Sethi, Vickson and Derzko \cite{BS} and Lasry and Lions \cite{LL}%
.

Our motivation in studying this equation comes from the recent work of \cite%
{DP}, where the author obtained non-positive radial solutions for the
equation (\ref{tpes}) and where we postulate an open problem regarding the
existence of positive solution for this equation. Another goal of this paper
is to improve the production planning model given in \cite{BS}, \cite{DP},
and to give a verification result, i.e., show that the solution of the
equation yields the optimal production.

To the best of our knowledge, the first mathematical result about the
existence of positive solution for the semilinear equation (\ref{tpes})
seems to be due to \cite{BS} for the case $N=1$ and no results on existence
of positive solutions are known for the case $N>1$. This should not surprise
us since there are some difficulties in analyzying this class of problems in 
$\mathbb{R}^{N}$, with $N\geq 1$, which will be revealed in the following
sections organized as follows. In Section \ref{mr}, we give our main theorem
regarding the existence of positive solution for the problem (\ref{tpes})
and its proof. The Section \ref{ppp} contains the application to a
production planning problem. Section \ref{avr}, presents a verification
result. In Section \ref{sc} we obtain a closed form solution for our
equation in a special case.

\section{Main Result \label{mr}}

The following Theorem is the main result of our paper.

\begin{theorem}
\label{msg}There exists a positive convex function $z\left( x\right) \in
C^{2}\left( \mathbb{R}^{N}\right) $ satisfying (\ref{tpes}). Moreover, the
following estimates hold%
\begin{eqnarray}
0 &\leq &z\left( x\right) \leq \left\vert x\right\vert ^{2}+1+\frac{N}{%
\alpha }\sigma ^{2},\text{ for }x\in \mathbb{R}^{N},  \label{ineq} \\
\left\vert \nabla {z}(x)\right\vert &\leq &C(1+\left\vert x\right\vert ),%
\text{ for }x\in \mathbb{R}^{N}\text{ and for some positive constant }C.
\label{sass}
\end{eqnarray}
\end{theorem}

To achieve our goal, which is establishing the above theorem, we prove the
following lemma:

\begin{lemma}
\label{ech}The partial differential equation with gradient term (\ref{tpes})
is equivalent to the semilinear elliptic equation%
\begin{equation}
-\Delta u\left( x\right) +\frac{1}{\sigma ^{4}}\left\vert x\right\vert
^{2}u\left( x\right) =-\frac{2\alpha }{\sigma ^{2}}u\left( x\right) \ln
u\left( x\right) \text{, for }x\in \mathbb{R}^{N}.  \label{tpe}
\end{equation}
\end{lemma}

\paragraph{\textbf{Proof.}}

With the change of variable $u\left( x\right) =e^{\frac{-z\left( x\right) }{%
2\sigma ^{2}}}$, the equation (\ref{tpe}) becomes (\ref{tpes}). Notice that
the above process is invertible, so (\ref{tpes}) and (\ref{tpe}) are
equivalent.

In the next theorem we prove the existence of a solution $u\left( x\right)
\in C^{2}\left( \mathbb{R}^{N}\right) $ for the problem (\ref{tpe}), such
that $0<u\left( x\right) \leq 1$ for all $x\in \mathbb{R}^{N}$.

\begin{theorem}
\label{mss}There exist functions $\underline{u},$ $\overline{u}:\mathbb{R}%
^{N}\rightarrow (0,1]$ of class $C^{2}(\mathbb{R}^{N})$ such that 
\begin{equation}
\left\{ 
\begin{array}{l}
-\Delta \underline{u}\left( x\right) +\frac{1}{\sigma ^{4}}\left\vert
x\right\vert ^{2}\underline{u}\left( x\right) \leq -\frac{2\alpha }{\sigma
^{2}}\underline{u}\left( x\right) \ln \underline{u}\left( x\right) \text{, }%
x\in \mathbb{R}^{N}, \\ 
-\Delta \overline{u}\left( x\right) +\frac{1}{\sigma ^{4}}\left\vert
x\right\vert ^{2}\overline{u}\left( x\right) \geq -\frac{2\alpha }{\sigma
^{2}}\overline{u}\left( x\right) \ln \overline{u}\left( x\right) ,\text{ }%
x\in \mathbb{R}^{N}, \\ 
\underline{u}\left( x\right) \leq \overline{u}\left( x\right) \text{ for }%
x\in \mathbb{R}^{N}.%
\end{array}%
\right.  \label{tpls}
\end{equation}%
Moreover, for such functions $\underline{u}$, $\overline{u}$ there exists a
function $u\left( x\right) \in C^{2}\left( \mathbb{R}^{N}\right) $ with $%
\underline{u}\left( x\right) \leq u\left( x\right) \leq \overline{u}\left(
x\right) $ in $\mathbb{R}^{N}$ and satisfying (\ref{tpe}).
\end{theorem}

Before giving the proof of the above theorem, we point that the function $%
\underline{u}$ (resp. $\overline{u}$) is called a sub-solution (resp.
super-solution) for the problem (\ref{tpe}).

\subparagraph{\textbf{Proof.}}

In the following we construct the functions $\underline{u},$ $\overline{u}.$
We adopt the idea of Bensoussan, Sethi, Vickson and Derzko \cite{BS}, for
the one dimensional case. More exactly, we find $a$, $b$ and $c$ such that $%
\underline{u}\left( x\right) =e^{a\left\vert x\right\vert ^{2}+b\left\vert
x\right\vert +c}$, is a sub-solution for the problem (\ref{tpe}). A simple
calculation shows that we can take $a=-\frac{1}{2\sigma ^{2}}$, $b=0$ and $%
c=-\frac{1}{2\sigma ^{2}}-\frac{N}{2\alpha }$, to provide the sub-solution
mentioned.

To construct a super-solution it is useful to remember that $\ln 1=0$ and
then a simple calculation shows that $\overline{u}\left( x\right) =1$, is a
super-solution of the problem (\ref{tpe}).

Until now, we constructed the corresponding sub- and super-solutions
employed in the one dimensional case by \cite{BS}. Clearly, (\ref{tpls})
holds and then in Theorem \ref{mss} it remains to prove that there exists $%
u\left( x\right) \in C^{2}\left( \mathbb{R}^{N}\right) $ with $\underline{u}%
\left( x\right) \leq u\left( x\right) \leq \overline{u}\left( x\right) $ in $%
\mathbb{R}^{N}$ satisfying (\ref{tpe}).

To do this, let $B_{k}=\left\{ x\in \mathbb{R}^{N}\left\vert \left\vert
x\right\vert <k\right. \right\} $ be the ball centered at the origin and of
radius $k=1,2,...$. We consider the problem%
\begin{equation}
\left\{ 
\begin{array}{l}
-\Delta u\left( x\right) +\frac{1}{\sigma ^{4}}\left\vert x\right\vert
^{2}u\left( x\right) =-\frac{2\alpha }{\sigma ^{2}}u\left( x\right) \ln
u\left( x\right) \text{, }x\in B_{k}, \\ 
u\left( x\right) =\underline{u}_{k}\left( x\right) \text{, }x\in \partial
B_{k},%
\end{array}%
\right.  \label{ball}
\end{equation}%
where $\underline{u}_{k}=\underline{u}_{\left\vert B_{k}\right. }$. In a
similar way, we define $\overline{u}_{k}=\overline{u}_{\left\vert
B_{k}\right. }$ and then $\underline{u}_{k}$, $\overline{u}_{k}\in
C^{2}\left( \overline{B}_{k}\right) $.

Observing that $\underset{x\in \mathbb{R}^{N}}{\inf }\underline{u}\left(
x\right) \leq \underset{x\in \overline{B}_{k}}{\min }\underline{u}_{k}\left(
x\right) $ and $\underset{x\in \mathbb{R}^{N}}{\sup }\overline{u}\left(
x\right) \geq \underset{x\in \overline{B}_{k}}{\max }\overline{u}_{k}\left(
x\right) $, a result of Kazdan and Kramer \cite{KK}, proves the existence of
a solution $u_{k}\in C^{2}\left( B_{k}\right) \cap C\left( \overline{B}%
_{k}\right) $ satisfying the problem (\ref{ball}). The function $u_{k}$ also
satisfies $\underline{u}_{k}\left( x\right) \leq u_{k}\left( x\right) \leq 
\overline{u}_{k}\left( x\right) $ for all $x\in \overline{B}_{k}$. By a
standard regularity argument based on Schauder estimates, see Noussair and
Swanson \cite[Lemma 3.2, p. 124]{NS} for details, we can see that for all
integers $k\geq n+1$ there are $\alpha \in \left( 0,1\right) $ and a
positive constant $C_{1}$ independent of $k$, such that 
\begin{equation}
u_{k}\in C^{2,\alpha }\left( \overline{B}_{n}\right) \text{ and }\left\vert
u_{k}\right\vert _{C^{2,\alpha }\left( \overline{B}_{n}\right) }<C_{1},
\label{b1}
\end{equation}%
where $\left\vert \circ \right\vert _{C^{2,\alpha }}$ is the usual norm of
the space $C^{2,\alpha }\left( \overline{B}_{n}\right) $. Moreover, there
exists a constant $C_{2}$ independent of $u_{k}$ such that%
\begin{equation}
\underset{x\in \overline{B}_{n}}{\max }\left\vert \nabla u_{k}\left(
x\right) \right\vert \leq C_{2}\underset{x\in \overline{B}_{k}}{\max }%
\left\vert u_{k}\left( x\right) \right\vert .  \label{b2}
\end{equation}%
The information from (\ref{b1}) and (\ref{b2}) implies that \{$\nabla u_{k}$%
\} as well as \{$u_{k}$\} are uniformly bounded on $\overline{B}_{n}$. Using
the compactness of the embedding $C^{2,\alpha }\left( \overline{B}%
_{n}\right) \hookrightarrow C^{2}\left( \overline{B}_{n}\right) $, enables
us to define the subsequence $u_{n}^{k}:=u_{k\left\vert B_{n}\right. },$ for
all $k\geq n+1$. Then for $n=1,2,3,...$ there exist a subsequence $%
\{u_{n}^{k_{nj}}\}_{k\geq n+1,j\geq 1}$ of $\{u_{n}^{k}\}_{k\geq n+1}$ and a
function $u_{n}$ such that $u_{n}^{k_{nj}}\rightarrow u_{n}$, uniformly in
the $C^{2}\left( \overline{B}_{n}\right) $ norm. More exactly, we get
through a well-known diagonal process that%
\begin{eqnarray*}
\mathbf{u}_{1}^{k_{11}}\text{, }u_{1}^{k_{12}}\text{, }u_{1}^{k_{13}}\text{, 
}... &\longrightarrow &u_{1}\text{ in }C^{2}\left( \overline{B}_{1}\right) ,
\\
u_{2}^{k_{21}}\text{, }\mathbf{u}_{2}^{k_{22}}\text{, }u_{2}^{k_{23}}\text{, 
}... &\longrightarrow &u_{2}\text{ in }C^{2}\left( \overline{B}_{2}\right) ,
\\
&&...
\end{eqnarray*}%
Since $\mathbb{R}^{N}=\underset{n=1}{\overset{\infty }{\cup }}B_{n}$, we can
define the function $u:\mathbb{R}^{N}\rightarrow \left[ 0,\infty \right) $
such that $u\left( x\right) =\lim_{n\rightarrow \infty }u_{n}\left( x\right) 
$. From the regularity theory the solution $u$ belongs to $C^{2}\left( 
\mathbb{R}^{N}\right) $, satisfies (\ref{tpe}) and the function $u$ also
satisfies $\underline{u}\left( x\right) \leq u\left( x\right) \leq \overline{%
u}\left( x\right) $, for all $x\in \mathbb{R}^{N}$.

\paragraph{\textbf{Proof of Theorem \protect\ref{msg}}}

The existence of solutions is proved by Lemma \ref{ech} and Theorem \ref{mss}%
. Then it remains to prove (\ref{ineq}), (\ref{ineq}) and that $z\left(
x\right) $ is a convex function.

A recapitulation, the change of variables say that $z\left( x\right)
=-2\sigma ^{2}\ln u\left( x\right) $, is a solution for (\ref{tpes}).
Observing that%
\begin{equation*}
\underline{u}\left( x\right) =e^{-\frac{\left\vert x\right\vert ^{2}}{%
2\sigma ^{2}}-\frac{1}{2\sigma ^{2}}-\frac{N}{2\alpha }}\leq u\left(
x\right) \leq \text{ }\overline{u}\left( x\right) =1\text{, }x\in \mathbb{R}%
^{N},
\end{equation*}%
it follows that%
\begin{equation}
0\leq z\left( x\right) \leq \left\vert x\right\vert ^{2}+1+\frac{N}{\alpha }%
\sigma ^{2}\text{, for }x\in \mathbb{R}^{N}.  \label{qua}
\end{equation}%
By the same arguments as in \cite[Theorem 3, p. 278]{AL} the solution $%
z\left( x\right) $ is convex. Since $z\left( x\right) $ verifies (\ref{qua})
the inequality (\ref{sass}) follows from \cite[Lemma 1, p. 24]{E} (see also
the arguments in \cite[Theorem 1, p. 236]{EP}). The proof is completed.

\section{Production Planning problem \label{ppp}}

As we already mentioned the equation studied is appearing in a stochastic
production planning problem. Indeed, let $\mathbb{R}^{N}$ ($N\geq 1$) be the 
$N$ dimensional Euclidean space and consider a factory producing $N$
homogeneous goods and having an inventory warehouse. Define the following
quantities:

1. $p\left( t\right) =\left( p_{1}(t),...,p_{N}(t)\right)$ represents the
production at time $t$ (control variable);

2. $p^{0}=\left( p_{1}^{0},...,p_{N}^{0}\right) $ represents factory optimal
production level;

3. $y\left( t\right) =\left( y_{1}(t),...,y_{N}(t)\right) $ denotes the
inventory level for production rate at time $t$ (state variable);

4. $l=\left( l_{1},...,l_{N}\right) $ denote the factory-optimal inventory
level;

5. $c$ represents production cost coefficient;

6. $h$ is the inventory holding cost coefficient;

7. $\xi =\left( \xi _{1},...,\xi _{N}\right) $ represents the constant
demand rate at time $t$;

8. $\sigma$ is a positive diffusion coefficient;

9. $\alpha >0$ is the constant discount rate;

10. $y_{i}^{0}$ is the initial inventory level;

11. $w=\left( w_{1},...,w_{N}\right) $ is a $N$-dimensional Brownian motion
on a complete probability space $(\Omega ,\mathcal{F},P)$ endowed with the
natural completed filtration $\{\mathcal{F}_{t}\}_{0\leq t\leq T}$, where $%
T=\infty$ is the length of planning period (we deal with the infinite
horizon case), and the filtration is generated by the standard Brownian
motion process $w.$

We now state the conditions of the model. The first condition is the state
dynamic equation for the inventory level stated as an It\^{o} stochastic
differential equation 
\begin{equation}
dy_{i}\left( t\right) =\left( p_{i}-\xi _{i}\right) dt+\sigma dw_{i}\text{, }%
y_{i}\left( 0\right) =y_{i}^{0}\text{, }i=1,...,N.  \label{sde}
\end{equation}%
The diffusion part of this equation is interpreted as "sales returns,"
"inventory spoilage," etc. which are random in nature.

The inventory production control problem is to choose the factory production
as to minimize the following cost functional

\begin{equation}
J\left( p_{1},...,p_{N}\right) :=\text{ }E\int_{0}^{\infty }(cf_{1}(\left(
p(t)-p^{0}\right) )+hf_{2}(y(t)-l))e^{-\alpha t}dt,  \label{ipcp}
\end{equation}%
where $f_{1}(x)=f_{2}\left( x\right) =\left\vert x\right\vert ^{2}$ is the
quadratic loss function.

To simplify the presentation we assume that $p^{0}\left( t\right) =l=\left(
0,...,0\right) $ and $h=c=1$. This assumption makes perfect sense if we
consider the deviations from the factory-optimal production level and
deviations from the factory-optimal inventory level. In light of this the
deviations are allowed to be negative. The aim is to minimize the stochastic
production planning problem 
\begin{equation}
\inf \{J\left( p_{1},...,p_{N}\right) \left\vert \,\,p_{i}\text{ }\forall
i=1,2,...,N\right. \}\text{, }  \label{pp}
\end{equation}%
with the inventory level subject to the It\^{o} equation (\ref{sde}).

Let $z\left( x\right) =z\left( x_{1},...,x_{N}\right) $ denote the expected
current-valued value of the control problem (\ref{sde})-(\ref{ipcp}) with
initial value $\left( x_{1},...,x_{N}\right) .$ In order to solve this
stochastic production planning problem we apply the martingale principle:
that is, we search for a function $U\left( x\right) $ such that the
stochastic process $M^{c}(t)$ defined below%
\begin{equation*}
M^{c}\left( t\right) =e^{-\alpha t}U\left( y\left( t\right) \right)
-\int_{0}^{t}[cf_{1}(p(s)-p^{0})+hf_{2}(y(s)-l)]e^{-\alpha s}ds,
\end{equation*}%
is supermartingale for all $p\left( t\right) =\left(
p_{1}(t),...,p_{N}(t)\right) $ and martingale for the optimal control $%
p^{\ast }\left( t\right) =\left( p_{1}^{\ast }(t),...,p_{N}^{\ast
}(t)\right) $. Then, it can be shown that $-U\left( x\right) =z\left(
x\right) $ is $C^{2}\left[ 0,\infty \right) $ and satisfies the
Hamilton-Jacobi-Bellman equation (HJB) formally associated with the problem (%
\ref{sde})-(\ref{pp}) 
\begin{equation}
\alpha z-\frac{\sigma ^{2}}{2}\Delta z+\xi \nabla z-\left\vert x\right\vert
^{2}=\inf \{p\nabla z+\left\vert p\right\vert ^{2}\left\vert p_{i}\text{ }%
\forall i=1,2,...,N\right. \},  \label{hjb}
\end{equation}%
where $z:=z\left( x_{1},...,x_{N}\right) $ is the corresponding value
function. The first order conditions yield the optimal candidate $p^{\ast
}\left( t\right) =\left( p_{1}^{\ast }(t),...,p_{N}^{\ast }(t)\right) $, by 
\begin{equation}
p_{i}^{\ast }(t)=-\frac{1}{2}\frac{\partial z}{\partial x_{i}}\left(
y_{1}^{\ast }(t),...,y_{N}^{\ast }(t)\right) \text{ for }i=1,...,n\text{.}
\label{optp}
\end{equation}%
and 
\begin{equation}
d{y_{i}^{\ast }}\left( t\right) =\left( {p_{i}^{\ast }}-\xi _{i}\right)
dt+\sigma dw_{i}\text{, }{y_{i}^{\ast }}\left( 0\right) ={y_{i}^{\ast 0}}%
\text{, }i=1,...,N.  \label{sde*}
\end{equation}

We point that the solution of (\ref{hjb}) equation is used to test
controller for optimality and equation (\ref{optp}) is used to construct a
feedback controller.

We consider the case $\xi =\left( 0,...,0\right) $ which makes sense if
deviation from the constant demand rate is taken into account. Then, this
equation (\ref{hjb}) can be simplified by noting that 
\begin{equation*}
\inf \{p\nabla z+\left\vert p\right\vert ^{2}\left\vert p_{i}\geq 0\text{ }%
\forall i=1,2,...,N\right. \}=-\frac{1}{4}\left\vert \nabla z\right\vert
^{2},
\end{equation*}%
so that equation (\ref{hjb}) can be written as%
\begin{equation}
-2\sigma ^{2}\Delta z+\left\vert \nabla z\right\vert ^{2}+4\alpha
z=4\left\vert x\right\vert ^{2}\text{ for }x\in \mathbb{R}^{N}.
\label{tpes1}
\end{equation}%
which is the same as equation (\ref{tpes}).

\section{A Verification Result\label{avr}}

In this section we establish the optimality of control $\left( p_{1}^{\ast
},...,p_{N}^{\ast }\right) $ given by (\ref{optp}) and (\ref{sde*}). The
verification theorem proceeds with the following steps:

\textbf{First Step:} The system of SDEs (\ref{sde*}) with $\left(
p_{1}^{\ast },...,p_{N}^{\ast }\right) $ given by (\ref{optp}) has a weak
solution via Girsanov Theorem in light of (\ref{sass}) (see Section 3.5 in 
\cite{KS} for more on this). Indeed this is true since in light of the
assumption the Novikov condition holds true on small intervals, and Girsanov
Theorem can be extended by an induction argument to arbitrarily large
intervals.

\textbf{Second Step:} Let $y_{t}^{\ast }$ be the inventory level
corresponding to $\left( p_{1}^{\ast },...,p_{N}^{\ast }\right) $ given by (%
\ref{optp}). In light of (\ref{sass}) one can get using the arguments
appearing in the proof of Theorem 5.2.1 from \cite{Oks}, we get for the
optimal control candidate%
\begin{equation}
E\left\vert y_{t}^{\ast }\right\vert ^{2}\leq C_{1}e^{C_{2}t},  \label{est}
\end{equation}%
for some positive constants $C_{1},{C_{2}}${.} Indeed, by integrating (\ref%
{sde*}), and by applying the expectation operator to $\left\vert y_{t}^{\ast
}\right\vert ^{2},$using Cauchy Schwarz inequality, employing (\ref{sass}),
and Gronwall inequality yields (\ref{est}).

\textbf{Third Step:} The set of acceptable production rates that we consider
is encompassing production rates and inventory levels for which%
\begin{equation*}
J\left( p_{1},...,p_{N}\right) :=\text{ }E\int_{0}^{\infty
}(f_{1}(p(t))+f_{2}(y(t)))e^{-\alpha t}dt<\infty ,
\end{equation*}%
and the following transversality condition $\lim_{t\rightarrow \infty
}Ee^{-\alpha t}\left\vert y_{t}\right\vert ^{2}=0$, is met. The set of
acceptable production rates is non empty. Because of (\ref{sass}) and
estimate (\ref{est}), the candidate optimal control $\left( p_{1}^{\ast
},...,p_{N}^{\ast }\right) $ verifies that $J\left( p_{1}^{\ast
},...,p_{N}^{\ast }\right) <\infty ,$ for $\alpha $ large enough. Moreover,
for $\alpha $ large enough the transversality condition is met because of (%
\ref{ineq}) and (\ref{est}). Also the control $p_{1}=0,...,p_{N}=0,$ is an
acceptable control. In light of the quadratic estimate on the value function
(see (\ref{ineq}) in the main theorem), the transversality condition implies
that $\lim_{t\rightarrow \infty }Ee^{-\alpha t}U(y_{t}^{\ast })=0$.

\textbf{Fourth Step:} Recall that $U(x)=-z(x),$ where $z$ is the solution of
(\ref{tpes}), and%
\begin{equation*}
M^{c}\left( t\right) =e^{-\alpha t}U\left( y\left( t\right) \right)
-\int_{0}^{t}(f_{1}(p(u))+f_{2}(y(u)))e^{-\alpha u}du.
\end{equation*}%
Therefore, the It\^{o}'s Lemma yields for the optimal control candidate%
\begin{equation*}
dM^{c}\left( s\right) =e^{-\alpha s}\sigma \nabla {z}(y^{\ast }(s))dw(s).
\end{equation*}%
Consequently $M^{c}(t)$ is a local martingale. Moreover,%
\begin{eqnarray*}
E\int_{0}^{t}e^{-2\alpha s}\sigma ^{2}\left\vert \nabla {z}(y^{\ast
}(s))\right\vert ^{2}ds &\leq &CE\int_{0}^{t}e^{-2\alpha s}\sigma
^{2}\left\vert y^{\ast }(s)\right\vert ^{2}ds+C_{3} \\
&\leq &CC_{1}\int_{0}^{t}e^{-2\alpha s}e^{C_{2}s}+C_{3}<C_{4},
\end{eqnarray*}%
for some positive constants $C,C_{1},C_{2},C_{3},C_{4},$ and large enough $%
\alpha .$ This in turn makes $M^{c}(t)$ a (true) martingale.

\textbf{Fifth Step: }This step establishes the optimality of $\left(
p_{1}^{\ast },...,p_{N}^{\ast }\right) $. The martingale/supermartingale
principle yields%
\begin{equation*}
Ee^{-\alpha t}U\left( y^{\ast }\left( t\right) \right)
-E\int_{0}^{t}(f_{1}(p^{\ast }(u))+f_{2}(y^{\ast }(u)))e^{-\alpha u}du=U(x)
\end{equation*}%
and%
\begin{equation*}
Ee^{-\alpha t}U\left( y\left( t\right) \right)
-E\int_{0}^{t}(f_{1}(p(u))+f_{2}(y(u)))e^{-\alpha u}du\leq U(x)
\end{equation*}%
By passing $t\rightarrow \infty $ and using transversality condition we get
the optimality of $\left( p_{1}^{\ast },...,p_{N}^{\ast }\right) .$

\section{Special Case \label{sc}}

In the following we manage to obtain a closed form solution for our equation
given a special discount $\alpha $. That is, assume $\alpha =N\sigma ^{2}$.
Then, two solutions for the problem (\ref{tpe}) are%
\begin{eqnarray}
u\left( \left\vert x\right\vert \right) &=&e^{m\left( \left\vert
x\right\vert ^{2}+1\right) }\text{, }m=\frac{1}{4\sigma ^{2}}\left( \alpha +%
\sqrt{\alpha ^{2}+4}\right) ,  \label{neg} \\
u\left( \left\vert x\right\vert \right) &=&e^{m\left( \left\vert
x\right\vert ^{2}+1\right) }\text{, }m=\frac{1}{4\sigma ^{2}}\left( \alpha -%
\sqrt{\alpha ^{2}+4}\right) .  \label{poz}
\end{eqnarray}%
Let us point out that (\ref{neg}) implies 
\begin{equation*}
z\left( x\right) =-2\sigma ^{2}m\left( \left\vert x\right\vert ^{2}+1\right)
<0\text{ for all }x\in \mathbb{R}^{N}\text{,}
\end{equation*}%
and then $z\left( x\right) $ is the negative solution obtained in the paper 
\cite{DP} and (\ref{poz}) implies that%
\begin{equation}
z\left( x\right) =-2\sigma ^{2}m\left( \left\vert x\right\vert ^{2}+1\right)
>0\text{ for all }x\in \mathbb{R}^{N}\text{,}  \label{poz1}
\end{equation}%
i.e. $z\left( x\right) $ is the positive solution obtained with the above
procedure. For the production planning problem we choose the positive
solution, i.e., the one given in (\ref{poz1}). Let us notice that $z\left(
x\right) $ given in (\ref{poz1}) satisfies the standing assumption (\ref%
{sass}), thus the verification holds true.

\textbf{Acknowledgments.} Traian A. Pirvu acknowledges that this work was
supported by NSERC grant 5-36700.


\begin{thebibliography}{99}
\bibitem{AK} R. Akella and P.R. Kumar, \textit{Optimal control of production
rate in a failure prone manufacturing system}, IEEE Trans Automat Contr. 31
(1986) 116-126.

\bibitem{A} O. Alvarez, \textit{A quasilinear elliptic equation in} $\mathbb{%
R}^{N}$, Proc. Roy. Soc. Edinburgh Sect. A. 126 (1996), 911-921.

\bibitem{AL} O. Alvarez, J.-M. Lasry and P.-L. Lions, \textit{Convex
viscosity solutions and state constraints}, J. Math. Pures Appl., 16 (1997)
265-288.

\bibitem{BS} A. Bensoussan, S.P. Sethi, R. Vickson, N. Derzko, \textit{%
Stochastic production planning with production constraints}, SIAM J. Control
Optim. 22 (1984) 920--935.

\bibitem{DP} D.-P. Covei, \textit{Symmetric solutions for an elliptic
partial differential equation that arises in stochastic production planning
with production constraints}, Appl. Math. Comput. 350 (2019) 190--197.

\bibitem{E} L. C. Evans, \textit{Weak Convergence Methods for Nonlinear
Partial Differential Equations}, CBMS Regional Conference Series in
Mathematics, 74 (1990) 1-82.

\bibitem{EP} L. C. Evans, R. F. Gariepy, \textit{Measure Theory and Fine
Properties of Functions}, Chapman and Hall/CRC Published April 14, 2015
Textbook - 313 Pages.

\bibitem{KK} J.L. Kazdan, R.J. Kramer, \textit{Invariant criteria for
existence of solutions to second-order quasilinear elliptic equations},
Comm. Pure Appl. Math. 31 (1978) 619--645.

\bibitem{KS} I. Karatzas, and S. E. Shreve, Brownian motion and stochastic
calculus, 2nd Ed., Springer-Verlag, New York, 1991.

\bibitem{LL} J.M. Lasry and P.L. Lions, \textit{Nonlinear elliptic equations
with singular boundary conditions and stochastic control with state
constraints}, Math. Ann. (1989) 283-583.

\bibitem{NS} E. S. Noussair and C. A. Swanson, \textit{Positive solutions of
quasilinear elliptic equations in Exterior Domains}, J. Math. Anal. Appl. 75
(1980) 121-133.

\bibitem{Oks} B. Oksendal, Stochastic Differential Equations,
Springer-Verlag, 5th Ed, 2000.
\end{thebibliography}
\end{document}